# A tale of two integrals:
# the Probability and Ahmed's integrals

___________________________

## 1) Introduction

One of the main rewards for the researcher in Mathematics is the feeling of surprise experienced when finding unexpected links between different topics or mathematical objects. It is such a link we chanced upon with the Ahmed integral.

The Ahmed integral appeared for the first time in 2001, in the Problems section of the *American Mathematical Monthly*. It has been named after its proposer, Zafar Ahmed, who requested the readers to evaluate this strange-looking integral:

$$A = \int_0^1 \frac{arctg\sqrt{2+x^2}}{(1+x^2)\sqrt{2+x^2}}\,dx$$

It was shown that $A = 5\pi^2/96$ (see proofs in [1]).

Since then several papers have dealt with alternative proofs of this result, and even generalisations (see for instance [2]).

Now what we intent to do is to describe another approach to Ahmed's integral, which will lead us to an unexpected link with the Probability integral, defined as

$$\int_0^{+\infty} e^{-x^2}\,dx = \tfrac{1}{2}\sqrt{\pi}$$

## 2) From the Probability integral to Ahmed's integral

In [3] we have explained how a new calculus method allows us to compute easily the square of the Probability integral, and the computation of its higher powers seemed to us the next natural move to continue the exploration of this method. Therefore to understand what follows the reader should previously get acquainted with the content of [3], and particularly its formulae *(F1)*, where $\alpha$ is a positive number and the function *f(x,y)* continuous on a subset of $\mathbf{R^2}$ containing $(0,\alpha) \times (0,\alpha)$ :

$$\int_0^\alpha \int_0^\alpha f(x,y)\,dxdy = \int_0^1 dx \int_0^\alpha \beta\,\{f(\beta,\beta x)+f(\beta x,\beta)\}\,d\beta \qquad (F1)$$

and *(F2)*

$$\int_0^\alpha \cdots\cdots \int_0^\alpha f(x_1,\cdots x_n)\,dx_1\cdots dx_n = \int_0^1 \cdots\cdots \int_0^1 dx_1\cdots dx_n \int_0^\alpha \beta^{n-1}\{\Sigma\}\,d\beta \qquad (F2)$$

where $\Sigma$ is the sum of all the functions $\Phi_p = f(\beta x_1,\cdots,\beta x_{p-1},\beta,\beta x_{p+1},\cdots,\beta x_n)$ with $p$ running from 1 to $n$.

Now let us consider the third power of an integral:

$$P = \left(\int_0^\alpha g(x)\,dx\right)^3$$

By the fundamental theorem of Calculus we have:

$$\frac{dP}{d\alpha} = 3g(\alpha)\left(\int_0^\alpha g(x)\,dx\right)^2$$

By the substitution $x \longrightarrow \alpha x$, and formula *(F1)*, we deduce that

$$\frac{dP}{d\alpha} = 3g(\alpha)\alpha^2\left(\int_0^1 g(\alpha x)\,dx\right)^2 = 6\int_0^1 dx\int_0^1 \alpha^2 \beta\, g(\alpha)\,g(\alpha\beta)\,g(\alpha\beta x)\,d\beta$$

Then, by integrating with regard to $\alpha$, we obtain:

$$P = \left(\int_0^\alpha g(x)\,dx\right)^3 = 6\int_0^1 dx\int_0^1 d\beta\int_0^\alpha \gamma^2 \beta\, g(\gamma)\,g(\gamma\beta)\,g(\gamma\beta x)\,d\gamma$$

By the same technique, this time with the help of formula *(F2)* with $n = 3$, we can prove that

$$\left(\int_0^\alpha g(x)\,dx\right)^4 = 4!\int_0^1 dx\int_0^1 d\beta\int_0^1 d\gamma\int_0^\alpha \delta^3\gamma^2\beta\, g(\delta)\,g(\delta\gamma)\,g(\delta\gamma\beta)\,g(\delta\gamma\beta x)\,d\delta \qquad (1)$$

Now, if in this formula we set $g(x) = e^{-x^2}$, we obtain

$$\left(\int_0^{+\infty} e^{-x^2}\,dx\right)^4 = \frac{\pi^2}{16} = 24\int_0^1 dx\int_0^1 d\beta\int_0^1 d\gamma\int_0^{+\infty} \delta^3\gamma^2\beta\, e^{-\delta^2(1+\gamma^2+\gamma^2\beta^2+\gamma^2\beta^2 x^2)}\,d\delta$$

Noticing that $2\delta\,d\delta = d(\delta^2)$, we can integrate with respect to $\delta$ and obtain *(2)*

$$\frac{\pi^2}{16} = 12 \int_0^1 dx \int_0^1 d\beta \int_0^1 d\gamma \ \frac{\gamma^2 \beta}{(1+\gamma^2+\gamma^2\beta^2+\gamma^2\beta^2 x^2)^2}$$

Since $2\beta\,d\beta = d(\beta^2)$, next we can integrate with regard to $\beta$ and obtain

$$\frac{\pi^2}{16} = 6 \int_0^1 dx \int_0^1 d\gamma \ \frac{1}{(1+\gamma^2)(1+x^2)} - 6 \int_0^1 dx \int_0^1 d\gamma \ \frac{1}{(1+x^2)(1+\gamma^2(2+x^2))}$$

The value of the first double integral above is obviously $(arctg1)^2 = \pi^2/16$. As for the second one, if we integrate with regard to $\gamma$, we find

$$\int_0^1 \frac{arctg\sqrt{2+x^2}}{(1+x^2)\sqrt{2+x^2}}\,dx$$

that is none other than Ahmed's Integral. As for its value, after having performed all the possible simplifications, we do find that it is $5\pi^2/96$, as expected.

### REFERENCES

**1. Definitely an Integral** *Am. Math. Month.* **109** (August-September 2002) pp 670-1

2. J.M. Borwein, D.H. Bailey and R. Girgensohn , *Experimentation in Mathematics : Computational Paths to Discovery* , AK Peters Ltd (March 2004), pp 17-20

3. Juan Pla, A footnote to the theory of double integrals, *Math. Gaz.* **94** n° 530 (July 2010) pp 262-9

*Juan PLA, 315 rue de Belleville, 75019 Paris (France)*